\documentclass[12pt,leqno]{article}

\newtheorem{example}{Example}
\newtheorem{theorem}{Theorem}

\newtheorem{corollary}{Corollary}
\newtheorem{remark}{Remark}
\newtheorem{lemma}{Lemma}
\newfont{\bb}{msbm10 at 12pt}
\newfont{\bbb}{msbm10 at 16pt}
\newfont{\got}{eufm10 at 12pt}

\def\C{\hbox{\bb C}}
\def\Cg{\hbox{\bbb C}}
\def\R{\hbox{\bb R}}
\def\S{\hbox{\bb S}}
\def\P{\hbox{\bb P}}

\def\Z{\hbox{\bb Z}}
\def\Q{\hbox{\bb Q}}
\def\s{\hbox{\got s}}
\def\o{\hbox{\got o}}
\def\u{\hbox{\got u}}
\def\g{\hbox{\got g}}
\def\l{\hbox{\got l}}

\begin{document}

\title{On a new construction of special Lagrangian immersions in complex Euclidean space}       

\author{Ildefonso Castro\thanks{Research partially 
supported by a MCYT and FEDER grant BFM2001-2967} \and Francisco Urbano$^{*}$}

\date{}

\maketitle

\section{Introduction}

 Special Lagrangian submanifolds of complex Euclidean space $\C^n$ have been studied widely over the last few years. These submanifolds are volume minimizing and, in particular, they are minimal submanifolds. When $n=2$, special Lagrangian surfaces of $\C^2$ are exactly complex surfaces with respect to another orthogonal complex structure on $\R^4\equiv\C^2$.
A very important problem here, and even a good starting point to study the not well understood problem about singularities of special Lagrangian $n$-folds in Calabi-Yau $n$-folds, is finding (non-trivial) examples of special Lagrangian submanifolds in $\C^n$, since those are locally modeled on singularities of these. 

In [2], R. Harvey and H.B. Lawson constructed the first examples in $\C^n$, where we point out the Lagrangian catenoid ([1, Remark 1], [3, Theorem A] and [5, Theorem 6.4]) to show a method of construction of special Lagrangian submanifolds of $\C^n$ using an (n-1)-dimensional oriented minimal Legendrian submanifold of $\S^{2n-1}$ and certain plane curves  (for a better understanding, see Proposition A in section 2). Making use of this method, we include in section 2 three new examples of special Lagrangian submanifolds which have not been exposed yet.

In [5], more examples of special Lagrangian submanifolds of $\C^n$ were constructed, all of them with cohomogeneity one, that is, the orbits of the symmetry group of the submanifold are of codimension one. Examples of special Lagrangian cones of cohomogeneity two were given in [3] and [5]. Recently, in [6] and [7], D. Joyce studied in depth special Lagrangian 3-submanifolds invariant under $U(1)$, where in general $U(n)$ denotes the unitary group.

In this paper we make a modest contribution to the construction of special Lagrangian submanifolds of $\C^n$ in the following way. We consider the group $SO(p+1)\times SO(q+1)$, with $p+q+2=n$, acting as a subgroup of the isometries group of $\C^n$ as follows: 
\[
(A,B)\in SO(p+1)\times SO(q+1)\longmapsto \left(\begin{array}{c|c}

A &  \\
\hline
&  B
\end{array}\right)\in U(n),
\]
and we analyze the special Lagrangian submanifolds of $\C^n$ invariant under the above action. When either $p=1,q=0$ or $p=0,q=1$, up to congruences, we are considering the $U(1)$-invariant special Lagrangian $3$-submanifolds studied by Joyce in [6] and [7]. The orbits of this action are isometric to the product of two spheres $\S^p\times\S^q$ embedded in $\C^n$ as a product of two totally geodesic Legendrian submanifolds: $$\S^p\times\S^q\hookrightarrow\S^{2p+1}\times\S^{2q+1}\subset\C^{p+1}\times\C^{q+1}\equiv\C^n.$$
When $p,q\geq 2$ the special Lagrangian $(SO(p+1)\times SO(q+1))$-invariant submanifolds can be described in terms of the above totally geodesic Legendrian submanifolds and certain Lagrangian surfaces of $\C^2$ whose definition can be found in Theorem 1.
Following a similar idea to the explained before, in Theorem 1 the above examples are generalized changing the totally geodesic Legendrian submanifolds by a pair of minimal Legendrian submanidolds in odd-dimensional spheres. These Lagrangian surfaces of $\C^2$ that appear in our construction will be studied in section 4 and can be considered as a generalization of the special Lagrangian surfaces of $\C^2$.  It is possible to describe them locally in terms of the solutions of the nonlinear Cauchy-Riemann equation:
\[
 f_x=g_y , \quad (q+1)^2(f^2+x^2)^{\frac{q}{q+1}}\, g_x=-(p+1)^2(g^2+y^2)^{\frac{p}{p+1}}\, f_y,
\]
where $f,g:U\subset \R^2 \rightarrow \R$ are differentiable functions (except at singular points) 
and $p$, $q$ are nonnegative integer numbers.
Therefore we expect that hard work involving the analytic techniques developed by Joyce in [6] and [7] would be successful in order to get existence, uniqueness and regularity of weak singular solutions of the above partial differential equations on strictly convex domains.

In Corollary 1, we construct in detail a 4-parameter family of such a class of Lagrangian surfaces of $\C^2$ from two couples of plane curves. Finally, in Corollary 2, we describe explicitly a 2-parameter subfamily of cylinders that generalize a family of complex curves of $\C^2$ studied in [4] by D.Hoffman and R.Osserman.

\section{Special Lagrangian submanifolds in $\Cg^n$}
Let $\C^n=\{(z_1,\dots,z_n)\,/\,z_i\in\C, \, i=1,\dots,n\}$ be the complex Euclidean space of dimension $n$ endowed with the bilinear product
\[
(z,w)=\sum_{i=1}^nz_i\bar{w}_i, \quad\forall z,w\in\C^n.
\]
Then $\langle,\rangle=\Re (,)$ is the Euclidean metric on $\C^n$ and $\omega=-\Im (,)$ the Kaehler two-form on $\C^n$, which is given by $\omega(v,w)=\langle Jv,w\rangle$, where $J$ is the complex structure on $\C^n$. The two-form $\omega$ is exact and it is given by $\omega=d\Gamma$, where $\Gamma$ is the Liouville $1$-form on $\C^n$ defined by
\[
2\Gamma(v)=\langle v,Jz\rangle,
\]
for all $v\in T_z\C^n$, $z\in\C^n$.

We consider the following $1$-parameter family of complex volume forms on $\C^n$:
\[
\Omega_{\theta}=e^{-i\theta}\cdot dz_1\wedge\dots\wedge dz_n,\quad \theta\in[0,2\pi[.
\]
Then $\Re (\Omega_{\theta})=\cos \theta\,\Re (\Omega_0)+\sin\theta\,\Im (\Omega_0)$ are calibrations on $\C^n$ called {\em special Lagrangian calibrations with phase $\theta$} (see [2]).

An immersion $\phi:M^n\rightarrow\C^n$ of a oriented manifold $M$ is called {\em special Lagrangian with phase $\theta$} if $\phi$ is calibrated with respect to the calibration $\Re (\Omega_{\theta})$, i.e. $\phi^*\Re (\Omega_{\theta})$ is the volume form on $M$. It is well-known (see [2, Corollary 1.11]) that {\em given an immersion $\phi:M^n\rightarrow\C^n$ of an orientable manifold $M$, then $M$ admits an orientation making $\phi$ a special Lagrangian immersion with phase $\theta$ if and only if $\phi$ is a Lagrangian immersion, i.e. $\phi^*\omega=0$ and $\phi^*\Im(\Omega_{\theta})=0$}. 
In this context, it is interesting to recall two important properties. The first one (see [2, Proposition 2.17]) is that {\em a Lagrangian immersion $\phi:M^n\rightarrow\C^n$ of an orientable manifold is minimal, i.e. the mean curvature vector of $\phi$ vanishes, if and only if $\phi$ is special Lagrangian with phase $\theta$, for some $\theta\in [0,2\pi[$}. The second one is that {\em if $\phi:M^n\rightarrow\C^n$ is a special Lagrangian immersion with phase $\theta$ and $A\in U(n)$, then $A\phi:M^n\rightarrow\C^n$ is a special Lagrangian immersion with phase $\theta+\alpha$, where $e^{i\alpha}=\det A$}.

The space of oriented Lagrangian $n$-planes in $\C^n$ is identified with the symmetric space $U(n)/SO(n)$, where
the determinant map $\det:U(n)/SO(n)$ $\rightarrow\S^1$ is well defined. If $\phi:M^n\rightarrow\C^n$ is a Lagrangian immersion of an oriented manifold and $\nu:M\rightarrow U(n)/SO(n)$ its Gauss map, then $\det\circ\nu:M\rightarrow\S^1$ is written by $\det\circ\nu=e^{i\beta}$, for a function $\beta:M\rightarrow\R/2\pi\Z$ called the {\em Lagrangian angle map} of $\phi$. This map $\beta$ satisfies that $J\nabla\beta=nH$, where $H$ is the mean curvature of $\phi$ and, in addition,  $e^{i\beta(p)}=(\Omega_0)_p(\phi_*T_pM)$ for all $p\in M$. Hence, {\em 
$\phi$ is a special Lagrangian immersion with phase $\theta$ if and only if $\beta(p)=\theta$, $\forall p\in M$}.

\vspace{0.5cm}

Among the known constructions of special Lagrangian submanifolds in $\C^n$, we are particularly interested in the following. If $\S^{2n-1}\subset\C^{n}$ denotes the unit sphere, an immersion $\psi:N^{n-1}\rightarrow\S^{2n-1}$ of an $(n-1)$-dimensional manifold $N$ is called {\em Legendrian} if $\psi^*\Gamma=0$. This implies that $\psi^*\omega=0$ and so $\psi$ is an isotropic immersion in $\C^n$. 
Legendrian submanifolds in $\S^{2n-1}$ are closely related with Lagrangian submanifolds in
the $(n-1)$-dimensional complex projective space $\C\P^{n-1}$. Concretely, if $\Pi:\S^{2n-1}\rightarrow\C\P^{n-1}$ is the Hopf fibration and $\psi:N^{n-1}\rightarrow\S^{2n-1}$ a Legendrian immersion, then $\Pi\circ\psi:N^{n-1}\rightarrow\C\P^{n-1}$ is a Lagrangian immersion. It is remarkable that the induced metrics on $N$ by $\psi$ and $\Pi\circ\psi$ are the same and that the mean curvature vector of $\psi$ is the horizontal lift of the mean curvature vector of $\Pi\circ\psi$. Moreover, it is interesting to mention that {\em locally} every Lagrangian submanifold of $\C\P^{n-1}$ is the projection by $\Pi$ of a Legendrian submanifold of $\S^{2n-1}$. 

\vspace{0.2cm}

{\sc Proposition A} ([1], [3], [5])

{\em Let $\gamma:I\rightarrow\C$ be a regular curve and $\psi:N^{n-1}\rightarrow\S^{2n-1}$ a Legendrian immersion of an orientable manifold $N$. Then $\Phi:I\times N\rightarrow\C^n$, defined by
\[
\Phi(s,x)=\gamma(s)\psi (x),
\]
is a special Lagrangian immersion if and only if $\psi$ is minimal and $\gamma^n$ has curvature zero.}

\vspace{0.2cm}

 It happens that $\Phi$ has singularities in the points $(s,x)\in I\times N$ where $\gamma(s)=0$. Up to rotations in $\C^n$, the curve $\gamma$ can be taken in such a way that $\gamma^n(s)=c+is$ with $c\geq 0$ (i.e.  $\Re\,\gamma^n=c$). We denote the above curve $\gamma$ by $\gamma_c$. When $c=0$, $\gamma_0$ is the straight line $\Re z=0$ and the examples constructed in this way are cones with links $\psi$. When $c>0$, $\gamma_c$ can be parameterized by $\gamma_c(s)=\sqrt[2n]{s^2+c^2}\, e^{i\frac{\arctan s}{n \, c}}$ and the examples constructed in this way were given in [1, Remark 1], [3, Theorem A] and [5, Theorem 6.4].

If we take $\psi$ in Proposition A  as the totally geodesic Legendrian embedding of $\S^{n-1}$ into $\S^{2n-1}$, we obtain the special Lagrangian submanifold invariant under the diagonal action of $SO(n)$ on $\C^n$ (see Theorem 3.5 in [2]). If we take $\psi$ as the standard Legendrian embedding of an $(n-1)$-torus in $\S^{2n-1}$, then we obtain the examples given by Harvey and Lawson in [2, Theorem 3.1]. But there are other important special Lagrangian submanifolds coming from Proposition A which have not been exposed yet. We are going to describe some of them.

\begin{example} Let $\g\l(n,\C)$, $n\geq 2$, the vector space of complex $n$-matrixes and $\langle \langle A,B\rangle\rangle=\Re\,(\hbox{\rm Trace}\,A\bar{B}^t)$ the Euclidean metric on $\g\l(n,\C)$. For each real number $a\geq 0$,
\[
M_a=\{B\in\g\l(n,\C)\,/\,B\bar{B}^t=|\hbox{\rm det}\,B|^{2/n}I \, , \, \Re\,((\hbox{\rm det}\,B)^n)=a\}
\]
is a special Lagrangian submanifold with phase $\pi/2$ of the complex Euclidean space $\g\l(n,\C)$. Moreover, $M_a$ is the only special Lagrangian submanifold of $\g\l(n,\C)$ invariant under the action of the special unitary group $SU(n)$ on $\g\l(n,\C)$ given by
\begin{eqnarray*}
SU(n)\times\g\l(n,\C)&\longrightarrow&\g\l(n,\C)\\
(A,B)&\mapsto&AB.
\end{eqnarray*}
\end{example}

{\it Proof:\/} Suppose that $\Phi:M\rightarrow\g\l(n,\C)$ is an special Lagrangian immersion $SU(n)$-invariant. Then for each $B\in M$ and $X\in \s\u(n)$, ($\s\u(n)$ stands for the Lie algebra of $SU(n)$), the curve $t\mapsto e^{tX}B$ is contained in the submanifold $M$. So its tangent vector at $t=0$ must be tangent to $M$ at $B$, i.e. $XB\in T_BM$, for any $X\in \s\u(n)$. As $\Phi$ is Lagrangian, for any $X,Y\in\s\u(n)$, we have that
\[
0=\langle\langle XB,JYB\rangle\rangle =-\Im \,(\hbox{\rm Trace}\,XB\bar{B}^tY).
\]
Since $X$ and $Y$ are arbitrary fields in the Lie algebra of $SU(n)$, from the above equation it is easy to conclude that
\( B\bar{B}^t=\lambda I\), with $\lambda$ a positive real number,
and hence there exists a non-null complex number $z$ and a matrix $A\in SU(n)$ such that $B=zA$. As the dimension of $M$ is $n^2$, we can parameterize locally $\Phi$ by
\begin{eqnarray*}
\Phi:I\times SU(n)&\longrightarrow&\g\l(n,\C)\\
(t,A)&\mapsto&\frac{\gamma(t)}{\sqrt n}A,
\end{eqnarray*}
where $I$ is an open interval in $\R$ and $\gamma:I\rightarrow\C$ is a curve. Now, we use that
\begin{eqnarray*}
\psi:SU(n)&\longrightarrow&\S^{2n^2-1}\subset\g\l(n,\C)\\
A&\mapsto&\frac{1}{\sqrt n}A,
\end{eqnarray*}
is a minimal Legendrian embedding (see [8]), and so $\Phi$ is one of the examples given in Proposition A. Therefore $\gamma=\gamma_c$ for some $c\geq 0$, where $(\gamma_c)^{n^2}(s)=c+is$. Now it is easy to check that $\Phi(I\times SU(n))\subset M_a$, with $a=c/n^{n^2/2}$. This finishes the proof.$_{\diamondsuit}$
\begin{remark}
{\rm When $n=2$, $SU(2)=\S^3$, $\hbox{dim}\,M_a=4$ and the examples $M_a$ are congruent to the given by Harvey and Lawson in [2, Theorem 3.5].}
\end{remark}

\begin{example} Let $S(n,\C)=\{B\in\g\l(n,\C)\,/\,B=B^t\}$, $n\geq 2$, the vector space of symmetric complex $n$-matrixes and $\langle\langle \,,\,\rangle\rangle$ the Euclidean metric induced on $S(n,\C)$. For each real number $b\geq 0$,
\[
M_b=\{B\in S(n,\C)\,/\,B\bar{B}=|\hbox{\rm det}\,B|^{2/n}I, \, \Re\,((\hbox{\rm det}\,B)^n)=b\}
\]
is a special Lagrangian submanifold with phase $\pi/2$ of the complex Euclidean space $S(n,\C)$. Moreover, $M_b$ is the only special Lagrangian submanifold of $S(n,\C)$ invariant under the action of $SU(n)$ on $S(n,\C)$ given by
\begin{eqnarray*}
SU(n)\times S(n,\C)&\longrightarrow&S(n,\C)\\
(A,B)&\mapsto&ABA^t.
\end{eqnarray*}
\end{example}

The proof of this result is similar to the proof of Example 1 using in this case the minimal Legendrian embedding \begin{eqnarray*}
\psi:SU(n)/SO(n)&\longrightarrow&\S^{n^2+n-1}\subset S(n,\C)\\
A\cdot SO(n) &\mapsto&\frac{1}{\sqrt n}AA^t ,
\end{eqnarray*}
given also by Naitoh in [8].

\begin{remark}
{\rm When $n=2$, $SU(2)/SO(2)=\S^2$, $\hbox{dim}\,M_b=3$ and the examples $M_b$ are congruent to the given by Harvey and Lawson in [2, Theorem 3.5]. In addition, the special Lagrangian submanifolds given in Example 2 are the intersections with the complex subspace $S(n,\C)$ of $\g\l(n,\C)$ of the special Lagrangian submanifolds given in Example 1.}
\end{remark}
\begin{example} Let $\s\o(2n,\C)=\{B\in\g\l(2n,\C)\,/\,B+B^t=0\}$, $n\geq 1$, the vector space of skew-symmetric complex $2n$-matrixes and $\langle\langle\,,\,\rangle\rangle$ the Euclidean metric induced on $\s\o(2n,\C)$. For each real number $c\geq 0$,
\[
M_c=\{B\in\s\o(2n,\C)\,/\,B\bar{B}=-|\hbox{\rm det}\,B|^{1/n}I, \, \Re\,((\hbox{\rm det}\,B)^{2n})=c\}
\]
is a special Lagrangian submanifold with phase $\pi/2$ of the complex Euclidean space $\s\o(2n,\C)$. Moreover, $M_c$ is the only special Lagrangian submanifold of $\s\o(2n,\C)$ invariant under the action of $SU(2n)$ on $\s\o(2n,\C)$ given by
\begin{eqnarray*}
SU(2n)\times\s\o(2n,\C)&\longrightarrow&\s\o(2n,\C)\\
(A,B)&\mapsto&ABA^t.
\end{eqnarray*}
\end{example}
Again the proof is similar to the above ones, considering in this last case the minimal Legendrian embedding (cf. [8])
given by
\begin{eqnarray*}
\psi:SU(2n)/Sp(n)&\longrightarrow&\S^{4n^2-1}\subset\s\o(2n,\C)\\
A\cdot Sp(n)&\mapsto&\frac{1}{\sqrt {2n}}AJ_nA^t,
\end{eqnarray*}
where  
$J_n=\left[\left(\begin{array}{c|c}
0 & -I_n \\
\hline
I_n & 0 
\end{array}\right)\right]$.
\begin{remark}
{\rm 
When $n=1$, $SU(2)/Sp(1)=\S^1$, $\hbox{dim}\,M_c=2$ and the examples $M_c$ are once again congruent to the given by Harvey and Lawson in [2, Theorem 3.5]. The special Lagrangian submanifolds given in Example 3 are also the intersections with the complex subspace $\s\o(2n,\C)$ of $\g\l(2n,\C)$ of the special Lagrangian submanifolds given in Example 1.}
\end{remark}

\section{New examples of special Lagrangian submanifolds in $\Cg^n$}

In the next result we generalize the construction exposed in Proposition A in the following way. 
\begin{theorem}
Let $\phi=(\phi_1,\phi_2):\Sigma\rightarrow\C^2$ be a Lagrangian immersion of an orientable surface $\Sigma$ and  $\psi:M^{p}\rightarrow\S^{2p+1}\subset\C^{p+1}$, $\varphi:N^q\rightarrow\S^{2q+1}\subset\C^{q+1}$ Legendrian immersions   of orientable manifolds $M$ and $N$. Then the Lagrangian immersion 
\[ 
\Phi:\Sigma\times M\times N \longrightarrow \C^n=\C^{p+1}\times\C^{q+1}, \, n=p+q+2, \] 
defined by
\[
\Phi(u,x,y)=(\phi_1(u)\psi (x),\phi_2(u)\varphi(y))
\]
is a special Lagrangian immersion if and only if $\psi:M^p\rightarrow\S^{2p+1}$ and $\varphi:N^q\rightarrow\S^{2q+1}$ are minimal immersions and the Lagrangian angle map $\beta_{\phi}:\Sigma\rightarrow\R/2\pi\Z$ of $\phi$ satisfies
\[
\beta_{\phi}+p\, \arg \phi_1 +q\, \arg\phi_2={\rm constant},
\]
where $\arg \phi_i:\Sigma\rightarrow \R/2\pi\Z$ are defined by $\arg\phi_i(u)=$argument$\,(\phi_i(u))$, $\forall u\in \Sigma$, $i=1,2$.

Moreover, any special Lagrangian immersion in $\C^n$ invariant under the action of $SO(p+1)\times SO(q+1)$, with $p+q+2=n$ and $p,q\geq 2$, is congruent to an open subset of one of the above special Lagrangian submanifolds when $\psi$ and $\varphi$ are the totally geodesic Legendrian embeddings of $\,\S^p$ and $\,\S^q$ in $\,\S^{2p+1}$ and $\,\S^{2q+1}$ respectively.
\end{theorem}
These special Lagrangian immersions described in Theorem 1 have singularities in the points $(u,x,y)\in \Sigma\times M\times N$ where either $\phi_1(u)=0$ or $\phi_2(u)=0$.
\begin{remark}
{\rm In the limit case $n=2$, which implies $p=q=0$ and $\psi=\varphi\equiv 1\in\S^1$, the above equation means that the Lagrangian immersion $\phi:\Sigma\rightarrow \C^2$ is special Lagrangian. 

If $n=3$, then either $p=1$ and $q=0$ or $p=0$ and $q=1$. In both cases, $\psi$ or $\varphi$ must be a minimal Legendrian curve in $\S^3$, i.e. a great circle in $\S^3$, which can be parameterized, for example,  by $\psi(t)=\frac{1}{\sqrt 2}(e^{it},e^{-it})$. Hence the special Lagrangian immersions $\Phi's$ of Theorem 1 when $n=3$ are given, up to congruences, by
\begin{eqnarray*}
\Phi:\Sigma\times\R&\longrightarrow&\C^3\\
(u,t)&\mapsto&\left(\frac{\phi_1(u)}{\sqrt 2}e^{it},\frac{\phi_1(u)}{\sqrt 2}e^{-it},\phi_2(u)\right),
\end{eqnarray*}
with $\beta_{\phi}+\arg \phi_1=$ constant. These examples are invariant under the action of $U(1)\equiv SO(2)$ on $\C^3$ given by
\[
e^{it}\cdot(z_1,z_2,z_3)=(e^{it}z_1,e^{-it}z_2,z_3)
\]
and they have been studied in depth by Joyce in [6] and [7]. 
}
\end{remark}
\begin{remark}
{\rm The product of two Lagrangian immersions given in Proposition A  provide examples of special Lagrangian submanifolds of Theorem 1 whose corresponding Lagrangian surfaces are written as  $(t,s)\in \R^2\mapsto(\gamma(t),\eta(s))\in\C^2$, where $\gamma^{p+1} $ and $\eta^{q+1}$ are plane curves with curvature zero.}
\end{remark}
{\it Proof:\/} We start proving the first part of the result.

If $g_{\phi}$, $g_{\psi}$ and  $g_{\varphi}$  are the induced metrics on $\Sigma$, $M$ and $N$ by $\phi$, $\psi$ and  $\varphi$ respectively, the induced metric on $\Sigma\times M\times N$ by $\Phi$ is given by $g_{\phi} + |\phi_1|^2g_{\psi}+|\phi_2|^2g_{\varphi}$. Hence the singularities of the induced metric are in the points where either $\phi_1=0$ or $\phi_2=0$. Since $\phi$ is a Lagrangian immersion and $\psi$ and $\varphi$ are Legendrian ones, it is clear that the immersion $\Phi$ is Lagrangian. 

In order to determine when $\Phi$ is a special Lagrangian immersion, we are going to compute the Lagrangian angle map $\beta_{\Phi}$. Let $\{e_1,e_2\}$, $\{v_1,\dots,v_{p}\}$ and $\{w_1,\dots,w_{q}\}$  be oriented local orthonormal references on $\Sigma$, $M$ and $N$ respectively. Then 
\begin{eqnarray*}
u_1&=&(e_1,0,0),\, u_2=(e_2,0,0),\\
u_i&=&(0,|\phi_1|^{-1}v_{i-2},0),\, i=3,\dots,p+2,\\
u_j&=&(0,0,|\phi_2|^{-1}w_{j-p-2}),\, j=p+3,\dots,n,
\end{eqnarray*}
is an oriented local orthonormal reference on $\Sigma\times M\times N$. 
So
\[
\begin{array}{c}
e^{i\beta_{\Phi}}=\det\,\{\Phi_*(u_1),\dots,\Phi_*(u_n)\}= \frac{\phi_1^{p}\phi_2^{q}}{|\phi_1|^{p}|\phi_2|^{q}} \times \\
\\
\det\,\{\Phi_*(u_1),\Phi_*(u_2),
(\psi_{*}(v_1),0),.,(\psi_{*}(v_{p}),0),(0,\varphi_{*}(w_1)),.,(0,\varphi_{*}(w_{q}))\}.
\end{array}
\]
But, for $i=1,2 $:
\[
\Phi_*(u_i)=\phi_{1_*}(e_i)(\psi,0)+\phi_{2_*}(e_i)(0,\varphi).
\]
In this way we obtain that 
\[
\begin{array}{c}
e^{i\beta_{\Phi}}=e^{i(\beta_{\phi}+p\arg\phi_1+q\arg\phi_2)}\times \\ \\
\det\{(\psi,0),(0,\varphi),
(\psi_{*}(v_1),0),.,(\psi_{*}(v_{p}),0),(0,\varphi_{*}(w_1)),.,(0,\varphi_{*}(w_{q})\}.
\end{array}
\]
Hence we finally arrive at 
\[
e^{i\beta_{\Phi}(u,x,y)}=(-1)^{p}\, e^{i(\beta_{\phi}+p\arg\phi_1+q\arg\phi_2)(u)}\,\det B(x)\,\det C(y),
\]
where $B(x)$ and $C(y)$ are the matrixes
\[
B(x)=\{\psi,\psi_{*}(v_1),\dots,\psi_{*}(v_{p})\}(x)
\]
and
\[
C(y)=\{\varphi,\varphi_{*}(w_1),\dots,\varphi_{*}(w_{q})\}(y).
\]
Thus $\Phi$ is a special Lagrangian immersion, i.e. $\beta_{\Phi}$ is constant, if and only if $\beta_{\phi}+p\,\arg\phi_1+q\,\arg\phi_2$ is constant on $\Sigma$ and $\det B$ and $\det C$ are constant on $M$ and $N$
respectively. 

Taking now into account that $B$ is a unitary matrix because of the Legendrian character of $\psi$, if $v$ is a tangent vector to $M$ we get
\[
v(\det B)=\det B\, \hbox{Trace}\,(v(B)\bar{B}^t).
\]
So $\det B$ is constant on $M$ if and only if $\hbox{Trace}\,(v(B)\bar{B}^t)=0$ for any $v$. If  $\sigma_\psi$ is the second fundamental form of $\psi:M\rightarrow\S^{2p+1}$, we have that
\[
v(B)=\{\psi_{*}(v),\sigma_\psi(v,v_1)-\langle v,v_1\rangle\psi,\dots,\sigma_\psi(v,v_{p})-\langle v,v_{p}\rangle\psi\},
\]
and, using again that $\psi$ is a Legendrian immersion, we easily obtain 
\[
\hbox{Trace}\,(v(B)\bar{B}^t)=i \,p\langle H_\psi,Jv\rangle,
\]
where $H_\psi$ is the mean curvature of $\psi$. As a summary, we have proved that $\det B$ is constant if and only if $\psi$ is a minimal immersion. An analogous result can be also obtained for $C$ and $\varphi$. 

On the other hand, let $\Phi:M^n\rightarrow\C^n$ be a special Lagrangian immersion of an orientable manifold $M$ invariant under the action of $SO(p+1)\times SO(q+1)$. Let $p$ be any point of $M$ and let $z=(z_1,\dots,z_n)=\Phi(p)$. As $\Phi$ is invariant under the action of $SO(p+1)\times SO(q+1)$, for any matrix $X=(X_1,X_2)$ in the Lie algebra of $SO(p+1)\times SO(q+1)$, the curve $s\mapsto ze^{s\hat X}$ with
\[
\hat X=\left( \begin{array}{c|c}
\mbox{X}_1 &  \\
\hline
& \mbox{X}_2
\end{array}\right)
\]
lies in the submanifold. Thus its tangent vector at $s=0$ satisfies
\[
z\hat X\in\Phi_*(T_pM).
\]
Since $\Phi$ is a Lagrangian immersion, this implies that
\[
\Im (z\hat X\hat Y\bar{z}^t)=0
\]
for any matrices $X=(X_1,X_2)$, $Y=(Y_1,Y_2)$ in the Lie algebra of $SO(p+1)\times SO(q+1)$. 
As $p+1\geq 3$ and $q+1\geq 3$,  it is easy to see from the last equation that $\Re (z_1,\dots,z_{p+1})$ and $\Im (z_{1},\dots,z_{p+1})$ (respectively $\Re (z_{p+2},\dots,z_{n})$ and $\Im (z_{p+2},\dots,z_{n})$) are linear dependent. As $SO(p+1)$ acts transitively on $\S^p$ and $SO(q+1)$ acts transitively on $\S^q$, we obtain that $z$ is in the orbit (under the action of $SO(p+1)\times SO(q+1)$ described above) of the point $(z^0_1,0,\dots,0,z^0_{p+2},0,\dots,0)$, with $|z^0_1|^2=\sum_{i=1}^{p+1}|z_i|^2$ and $|z^0_{p+2}|^2=\sum_{j=p+2}^{n}|z_j|^2$. This implies that locally $\Phi$ is the orbit under the action of $SO(p+1)\times SO(q+1)$ of a surface in $\C^2\equiv \C^n\cap\{z_2=\dots=z_{p+1}=z_{p+3}=\dots=z_n=0\}$. Therefore $M$ is locally $\Sigma\times\S^p\times\S^q$, with $\Sigma$ a surface. Moreover, $\Phi$ is given by
\[
\Phi(u,x,y)=(\phi_1(u)x,\phi_2(u)y),
\]
where $\phi=(\phi_1,\phi_2)$ is a Lagrangian surface in $\C^2$. Finally, as $\Phi$ is a special Lagrangian surface, the result follows using the first part of this Theorem. $_\diamondsuit$

\section{The special class of Lagrangian surfaces }

As a consequence of Theorem 1, our interest in this section will be to study Lagrangian surfaces $\phi:\Sigma\rightarrow\C^2$ satisfying $\beta_{\phi}+p\arg\phi_1+q\arg\phi_2=$ constant, where $p,q$ are nonnegative integer numbers. Up to congruences, we can essentially  look for Lagrangian immersions
$\phi=(\phi_1,\phi_2):\Sigma\rightarrow\C^2$ satisfying $\beta_{\phi}+p\arg\phi_1+q\arg\phi_2=0$.

\subsection{An analytic scheme}

In a first attempt, we could follow some ideas of Joyce in [6] and [7]. 
We define $\hat{\phi}:\Sigma\rightarrow\C^2$ by
\[
\hat{\phi}=(\phi_1^{p+1},\phi_2^{q+1}).
\]
Then $\hat{\phi}$ is an immersion with singularities and it is easy to check that
\[
\Omega_0({\hat{\phi}_*T\Sigma})=e^{i(\beta_{\phi}+p\arg\phi_1+q\arg\phi_2)}.
\]
So, the condition $\beta_{\phi}+p\arg\phi_1+q\arg\phi_2=0$ implies that $\hat{\phi}^*(\Im \,\Omega_0)=0$. But $\phi$ to be Lagrangian has not a good translation to $\hat{\phi}$. We write locally $\hat{\phi}$ as the graph $\hat{\phi}(x,y)=(g(x,y)+iy,f(x,y)-ix)$ for certain functions $f,g$. Under these conditions, it is easy to check that $\phi$ is a Lagrangian immersion satisfying $\beta_{\phi}+p\arg\phi_1+q\arg\phi_2=0$ if and only if $f,g$ satisfy the following differential equations: 
\begin{eqnarray}
 f_x&=& g_y\\\nonumber
(q+1)^2(f^2+x^2)^{\frac{q}{q+1}}\, g_x&=&-(p+1)^2(g^2+y^2)^{\frac{p}{p+1}}\, f_y,
\end{eqnarray}
which can be considered as a generalization of the Cauchy-Riemann equations that corresponds to $p=q=0$.
The first equation just means that $\hat{\phi}^*(\Im \,\Omega_0)=0$ and the second one that $\phi$ is a Lagrangian immersion. The solutions $f(x,y)=ax+b$, $g(x,y)=ay+c$, $a,b,c\in\R$, of (1) provide the examples described in Remark 5.

The first equation allows to rewrite (1) in terms of a potential function $h$ such that $h_x=g$ and $h_y=f$. So (1) becomes in the following second order partial differential equation:
\begin{equation}
(q+1)^2(h_y^2+x^2)^{\frac{q}{q+1}}\, h_{xx}+(p+1)^2(h_x^2+y^2)^{\frac{p}{p+1}}\,h_{yy}=0.
\end{equation}
When $p=1$ and $q=0$, the equations (1) and (2) correspond to the equations (12) and (18)---with $a=0$--- in [7].
In our case, at points $(x,0)$ with $g(x,0)=0$ and at points $(0,y)$ with $f(0,y)=0$, the equation (1) is no longer elliptic
as well as in [7]. But if we slightly modify (1), and consequently (2), in order to avoid singular points, following [6]
we would be able to prove the analogous result to Theorem 7.9 in [6],
that is, {\em the Dirichlet problem for the second-order quasilinear elliptic equation 
\[
(q+1)^2(h_y^2+x^2+a_1^2)^{\frac{q}{q+1}}\, h_{xx}+(p+1)^2(h_x^2+y^2+a_2^2)^{\frac{p}{p+1}}\,h_{yy}=0,
\]
($a_1\neq0$, $a_2\neq 0$) is uniquely solvable in a strictly convex domain of $\R^2$.} Notice that when $a_1=a_2=0$ the above equation is nothing but (2). To get such a result for (2), in the spirit of Theorem 8.17 in [7], carries much more difficulty. The key point is centered in making detailed analytic estimates of $f$, $g$ and their derivatives in order to prove a priori interior estimates for them when $f$, $g$ are bounded solutions of the above equation with $a_1\neq 0$ and $a_2\neq 0$
on a domain in $\R^2$.
So we are fully convinced that after a large hard work following similar technical analytic tools developed in [7], it could be proved existence and uniqueness of {\em weak} solutions to the Dirichlet problem for (2) on strictly convex domains.

\subsection{A 4-parameter family in terms of curves}

Now we are going to construct, for each pair of nonnegative integer numbers $p$ and $q$,
a $4$-parameter family of Lagrangian surfaces $\phi:\Sigma\rightarrow\C^2$ satisfying $\beta_{\phi}+p\arg\phi_1+q\arg\phi_2=0$. The construction involves two couples of planar curves whose definitions and properties are studied in the following results. 

\begin{lemma}
For any $p,q\geq 0$, let $\alpha_a=(\alpha_1,\alpha_2):I\subset \R \rightarrow \C^2$ be the only curve solution of 
\begin{equation}
\alpha'_j\bar{\alpha}_j=i\,\bar{\alpha}_1^{p+1}\bar{\alpha}_2^{q+1}, \  j=1,2,
\end{equation}
satisfying the real initial conditions $\alpha_a(0)=a=(a_1,a_2),\, a_1,a_2>0$.
Then:
\begin{enumerate}
\item $|\alpha_1|^2-|\alpha_2|^2=a_1^2-a_2^2$.
\item $\Re (\alpha_1^{p+1}\alpha_2^{q+1})=a_1^{p+1}a_2^{q+1}$, i.e. $\alpha_1^{p+1}\alpha_2^{q+1}$ is a straight line in $\C$.
\item For $j=1,2$, 
$\bar{\alpha}_j(t)=\alpha_j(-t),\,\forall t\in I$. 
\item The curves $\alpha_j,\,j=1,2$, are embedded and can be parameterized by $\alpha_j(s)=\rho_j(s)e^{i\theta_j(s)}$, where
\begin{eqnarray*}
\rho_j(s)&=&\sqrt{s^2+a_j^2},\\
 \theta_j(s)&=&\int_0^s\frac{a_1^{p+1}a_2^{q+1}x\,dx}{(x^2+a_j^2)\sqrt{(x^2+a_1^2)^{p+1}(x^2+a_2^2)^{q+1}-a_1^{2(p+1)}a_2^{2(q+1)}}}.
\end{eqnarray*}
\end{enumerate}
\end{lemma}
{\it Proof:\/} From (3) it follows that
\[
\frac{d}{dt}\left(|\alpha_1|^2-|\alpha_2|^2\right)=2\Re\,(\alpha_1'\bar{\alpha}_1-\alpha_2'\bar{\alpha}_2)=0,
\]
which proves 1. Using (3) again, we obtain that
\[
\frac{d}{dt}\left(\alpha_1^{p+1}\alpha_2^{q+1}\right)=i\,|\alpha_1|^{2p}|\alpha_2|^{2q}\left((p+1)|\alpha_2|^2+(q+1)|\alpha_1|^2\right),
\]
and this implies 2. In order to prove 3., it is enough to check that the curves $\beta_j(t)=\alpha_j(-t)$ and  $\delta_j(t)=\bar{\alpha}_j(t),\,j=1,2$, satisfy the same differential equation, since $\beta_j(0)=\delta_j(0)=a_j$, $j=1,2$.
 
Finally we prove 4. If we put
\[
\alpha_j(t)=\rho_j(t)e^{i\theta_j(t)},\,j=1,2,
\]
the equation (3) is equivalent to the equations
\begin{eqnarray}
\rho_j'\rho_j=\rho_1^{p+1}\rho_2^{q+1}\sin ((p+1)\theta_1+(q+1)\theta_2), \quad j=1,2,
\\
\rho_j^2\theta'_j=a_1^{p+1}a_2^{q+1},
\quad j=1,2.
\end{eqnarray}
The equation (5) says that the functions $\theta_j$, $j=1,2$, are strictly increasing. This means, using (4), that $t=0$ is the only critical point of $\rho_j$, $j=1,2$. 
In fact, it is not difficult to show that $t=0$ is a minimum of $\rho_j,\,j=1,2$. 
We know from 3. that $\rho_j(-t)=\rho_j(t)$ and $\theta_j(-t)=-\theta_j(t)$, $j=1,2$. 
So it is easy to get that $\alpha_j$, $j=1,2$, are embedded curves and
we can reparametrize them through the following common change of parameter:
\[
s=h(t)=\sqrt{\rho_j^2(t)-a_j^2},\quad j=1,2.
\]
Using (4), (5) and 2. we deduce  that
\[
\frac{d\theta_j}{ds}=\frac{a_1^{p+1}a_2^{q+1}s}{\rho_j^2\sqrt{\rho_1^{2(p+1)}\rho_2^{2(q+1)}-a_1^{2(p+1)}a_2^{2(q+1)}}}, \quad j=1,2,
\]
which finishes the proof of 4. $_{\diamondsuit}$

\begin{lemma}
For any $p,q\geq 0$, let $\gamma_b=(\gamma_1,\gamma_2):I'\subset \R \rightarrow \C^2$ be the only curve solution of
\begin{equation}
\gamma'_j\bar{\gamma}_j=(-1)^{j-1}\, i\, \bar{\gamma}_1^{p+1}\bar{\gamma}_2^{q+1}, \, j=1,2,
\end{equation}
satisfying the real initial conditions
$\gamma_b(0)=b=(b_1,b_2), \,  b_1, b_2>0$.
Then:
\begin{enumerate}
\item $|\gamma_1|^2+|\gamma_2|^2=b_1^2+b_2^2$.
\item $\Re (\gamma_1^{p+1}\gamma_2^{q+1})=b_1^{p+1}b_2^{q+1}$, i.e. $\gamma_1^{p+1}\gamma_2^{q+1}$ is a straight line in $\C$.
\item For $j=1,2$, $\bar{\gamma}_j(t)=\gamma_j(-t),\,\forall t\in I'$.
\item If $\, n^n \, b_1^{2(p+1)}b_2^{2(q+1)} <   (p+1)^{p+1}(q+1)^{q+1} (b_1^2+b_2^2)^{n}$, $n=p+q+2$, the functions $|\gamma_j|,\,j=1,2,$ are periodic with the same period $T=T(b)$. In this case, for such b's, $\gamma_b$ is a closed curve if and only if 
\[
\frac{b_1^{p+1}b_2^{q+1}}{2\pi}\left ( \int_0^{T}\frac{dt}{|\gamma_1|^2(t)},\int_0^{T}\frac{dt}{|\gamma_2|^2(t)}\right ) \in\Q^2 .
\]
\item If $b=\frac{1}{\sqrt{\lambda}} \left (\sqrt{p+1},\sqrt{q+1}\right )$, 
with $\lambda^{p+q}=(p+1)^{p}(q+1)^{q}$ (that corresponds with the equality in the inequality of 4.),
the curve $\gamma_b$ is given by
\[
\gamma_b(t)=\frac{1}{\sqrt{\lambda}} 
\left ( \sqrt{p+1}\,e^{i\sqrt{\frac{q+1}{p+1}}\,t},\sqrt{q+1}\,  
e^{-i\sqrt{\frac{p+1}{q+1}}\, t}\right ).
\]

\end{enumerate}
\end{lemma}
\begin{remark}
{\rm The curve given in Lemma 2.5. defines an embedding $\gamma_b:\S^1\rightarrow\C^2$ such that 
\[
\gamma_b(\S^1)\!=\! \{(z,w)\in\C^2\,/\,  \lambda z^{p+1}w^{q+1}= \sqrt{(p+1)(q+1)}\,,\,(q+1)|z|^2=(p+1)|w|^2 \}.
\]
}
\end{remark}

{\it Proof:\/} We omit the proofs of 1., 2. and 3. because they are very similar to the corresponding ones in Lemma 1.  The proof of 5. is a straightforward computation. 
To prove 4. we proceed as follows. 
If we write
\[
\gamma_j(t)=\varrho_j(t)e^{i\vartheta_j(t)},\,j=1,2,
\]
the equation (6) is equivalent to the equations
\begin{eqnarray}
\varrho_j'\varrho_j=(-1)^{j-1}\varrho_1^{p+1}\varrho_2^{q+1}\sin ((p+1)\vartheta_1+(q+1)\vartheta_2), \quad j=1,2,
\\
\varrho_j^2\vartheta'_j=(-1)^{j-1}b_1^{p+1}b_2^{q+1},
\quad j=1,2.
\end{eqnarray}
Using 1. we can put 
\[
\varrho_1=\sqrt{b_1^2+b_2^2}\, \cos \nu, \quad \varrho_2=\sqrt{b_1^2+b_2^2}\, \sin \nu,
\]
for a certain function $\nu=\nu(s)\in [0,\pi/2]$. Then (7) and 2. leads to 
\[
(b_1^2+b_2^2)^2 \sin^2 \nu \cos ^2 \nu \, (\nu')^2= (b_1^2+b_2^2)^n\cos^{2(p+1)} \nu \, \sin^{2(q+1)} \nu - b_1^{2(p+1)}b_2^{2(q+1)}.
\]
We analyze the set of critical points of $\nu$ from the above o.d.e. First we observe that $\nu(2s_0-s)=\nu(s)$ if $\nu'(s_0)=0$. Let $x=\cos^2 \nu(s_0) \in [0,1]$. Then $x$ must satisfy $(b_1^2+b_2^2)^{n}x^{p+1}  \,(1-x)^{q+1} -b_1^{2(p+1)}b_2^{2(q+1)}=0$. Studying the roots of this equation in $[0,1]$ we conclude that $\nu$ has exactly two critical 
values only if $n^n\, b_1^{2(p+1)}b_2^{2(q+1)}  < (p+1)^{p+1}(q+1)^{q+1}(b_1^2+b_2^2)^{n}$. In this case,
$\varrho_1$ and $\varrho_2$ are periodic functions with the same period, say $T$.
If $\gamma_b $ is a periodic curve, then its period must be an integer multiple of $T=T(b)$. Using (8)
it is not complicated to get that $\gamma_b$ is a closed curve if and only if $\vartheta_j(T)\in 2\pi \Q$, $j=1,2$, what proves 4.
$_{\diamondsuit}$

We already have the tools for the construction of a 4-parameter family of Lagrangian surfaces in $\C^2$ verifying the desired property.
\begin{corollary}
For each pair of integer numbers $p,q\geq 0$, 
let 
$\alpha_a=(\alpha_1,\alpha_2)$, $a=(a_1,a_2), \,  a_1, a_2>0$, and $\gamma_b=(\gamma_1,\gamma_2)$, $b=(b_1,b_2), \,  b_1, b_2>0$, be the curves of $\C^2$ 
introduced in Lemmas 1 and 2.
Then the map
\[
\phi:I\times I'\subset \R^2 \longrightarrow\C^2,
\]
defined by
\[
\phi(t,s)=(\alpha_1(t)\gamma_1(s),\alpha_2(t)\gamma_2(s)),
\]
is a Lagrangian immersion satisfying $\beta_{\phi}+p\arg\phi_1+q\arg\phi_2=0$.
\end{corollary}

{\it Proof:\/} Taking into account (3) and (6), it follows that
\[
\left(\frac{\partial\phi}{\partial t},\frac{\partial\phi}{\partial s}\right)=
\sum_{j=1}^2\alpha'_j\bar{\alpha}_j\gamma_j\bar{\gamma}'_j=\alpha'_1\bar{\alpha}_1\sum_{j=1}^2\gamma_j\bar{\gamma}'_j=0.
\]
Thus $\phi$ is a Lagrangian immersion. It is clear that
\[
e^{i(p\arg\phi_1+q\arg\phi_2)}=\frac{\alpha_1^{p}\alpha_2^{q}\gamma_1^{p}\gamma_2^{q}}{|\alpha_1|^{p}|\alpha_2|^{q}|\gamma_1|^{p}|\gamma_2|^{q}}.
\]
On the other hand, from section 2 it follows
\[
e^{i\beta_{\phi}(t,s)}=\Omega_0(\phi_* T_{(t,s)}(I\times I')).
\]
So, using again (3) and (6), we have that
\[
e^{i\beta_{\phi}}=\frac{\det\,\{\frac{\partial\phi}{\partial t},\frac{\partial\phi}{\partial s}\}}{|\frac{\partial\phi}{\partial t}\wedge\frac{\partial\phi}{\partial s}|}=\frac{(|\alpha_1|^2|\gamma_2|^2+|\alpha_2|^2|\gamma_1|^2)\bar{\alpha}_1^{p}\bar{\alpha}_2^{q}\bar{\gamma}_1^{p}\bar{\gamma}_2^{q}}{|\frac{\partial\phi}{\partial t}\wedge\frac{\partial\phi}{\partial s}|}=\frac{\bar{\alpha}_1^{p}\bar{\alpha}_2^{q}\bar{\gamma}_1^{p}\bar{\gamma}_2^{q}}{|\alpha_1|^{p}|\alpha_2|^{q}|\gamma_1|^{p}|\gamma_2|^{q}}.
\]
Finally, we get
\[
e^{i(\beta_\phi +p\arg\phi_1+q\arg\phi_2)}=1,
\]
which proves the result.$_{\diamondsuit}$

By combining Corollary 1 and Theorem 1 we conclude that from Legendrian minimal immersions $\psi $ and $\varphi $ in odd-dimensional spheres and the curves $\alpha_a=(\alpha_1,\alpha_2)$
and $\gamma_b=(\gamma_1,\gamma_2)$ of Lemmas 1 and 2, the immersions
\[
(t,s,x,y) \longmapsto (\alpha_1(t)\gamma_1(s)\psi(x),\alpha_2(t)\gamma_2(s)\varphi(y))
\]
provides a 4-parameter family of special Lagrangian submanifolds in $\C^n$.

\begin{corollary}
For each $a=(a_1,a_2)$ with $a_1,  a_2>0$ and each pair of integer numbers $p,q\geq 0$, 
\[
\left\{z\in\C^2\,/\,\frac{|z_1|^2}{p+1}-\frac{|z_2|^2}{q+1}=\frac{a_1^2-a_2^2}{\lambda},\,\Re\,(z_1^{p+1}z_2^{q+1})=\frac{ a_1^{p+1}a_2^{q+1}\sqrt{(p+1)(q+1)}}{\lambda}\right\}
\]
with $\lambda^{p+q}=(p+1)^{p}(q+1)^{q}$, 
is a Lagrangian surface $\Sigma_a$ in $\C^2$ satisfying $\beta_{\Sigma_a}+p\arg (\Sigma_a)_1+q\arg (\Sigma_a)_2=0$.
\end{corollary}
{\it Proof:\/} Taking in Corollary 1 any curve $\alpha_a$ from Lemma 1 and the curve given in Lemma 2,5., it is easy to prove this Corollary  making use of Remark 6.$_{\diamondsuit}$ 
\begin{remark}
{\rm The surfaces $\Sigma_a$ are topologically $\R\times\S^1$. Moreover, following Remark 2,  in the limit case $p=q=0$,
 $\Sigma_a$ defines a $2$-parameter family of special Lagrangian cylinders, which corresponds with the family of complex surfaces of $\C^2$ with finite total curvature $-4\pi$ given by Hoffman and Osserman in [4, Proposition 6.6, case 2]. }
\end{remark}

\vspace{1cm}

\begin{tabular}{ll}
{\sc addresses}: &  \\
(first author) & (second author) \\
Departamento de Matem\'{a}ticas & Departamento de Geometr\'{\i}a  \\
Escuela Polit\'{e}cnica Superior & y Topolog\'{\i}a \\
Universidad de Ja\'{e}n & Universidad de Granada\\
23071 Ja\'{e}n & 18071 Granada \\
SPAIN & SPAIN \\
 {\tt icastro@ujaen.es}& {\tt furbano@ugr.es}\\
\end{tabular}
\end{document}